\documentstyle[12pt]{article}
\newcommand{\qed}{\hfill\rule{4pt}{8pt}\par\vspace{\baselineskip}}

\newtheorem{de}{Definition}[section]
\newtheorem{lm}[de]{Lemma}

\newtheorem{pr}[de]{Proposition}
\newtheorem{co}[de]{Corollary}
\newtheorem{re}[de]{Remark}

\newtheorem{te}[de]{Theorem}
\newtheorem{ex}[de]{Example}

\def\Box{\mbox{$\sqcap\!\!\!\!\sqcup$}}

\def\ra{\rightarrow}
\def\lra{\longrightarrow}
\def\al{\alpha}

\def\bea{\begin{eqnarray*}}
\def\eea{\end{eqnarray*}}

\begin{document}
\title{Relative Regular Objects in Categories}
\author{ S.
D\u{a}sc\u{a}lescu$^1$, C. N\u{a}st\u{a}sescu$^1$,
A. Tudorache$^2$ and L. D\u{a}u\c{s}$^3$\\[2mm]
$^1$ University of Bucharest, Facultatea de Matematic\u{a},\\ Str.
Academiei 14, Bucharest 1, RO-010014, Romania\\
$^2$ Universitatea de \c{S}tiin\c{t}e Economice, Bucharest, Romania\\
$^3$ Universitatea de Construc\c{t}ii, Bucharest, Romania
 }

\date{}
\maketitle

\begin{abstract} We define the concept of a regular object with respect
to another object in an arbitrary category. We present basic
properties of regular objects and we study this concept in the
special cases of abelian categories and locally finitely generated
Grothendieck categories. Applications are given for categories of
comodules over a coalgebra and for categories of graded modules,
and a link to the theory of generalized inverses of matrices is
presented. Some of the techniques we use are new, since dealing
with arbitrary categories allows us to pass to the dual category.
\end{abstract}

\section{Introduction and preliminaries}

Von Neumann regular rings play a fundamental role in Ring Theory,
see \cite{go}. A ring $R$ is called von Neumann regular if for any
$a\in R$ there exists $b\in R$ such that $a=aba$. This concept was
generalized to modules in \cite{z}. A left module $M$ over the
ring $R$ is called regular if for any $m\in M$ there exists $g\in
Hom_R(M,R)$ such that $g(m)m=m$. Basic properties of regular
modules are developed in \cite{z}. Since a morphism $f\in
Hom_R(R,M)$ is uniquely given by an element $m\in M$, one can
reformulate the regular condition as follows. For any $f\in
Hom_R(R,M)$ there exists $g\in Hom_R(M,R)$ such that $f=f\circ
g\circ f$. This suggests the definition of a more general concept
of an $U$-regular object in a category, where $U$ is a given
object of the category. We study basic properties of regular
objects in categories, with special emphasis on abelian categories
and on locally finitely generated Grothendieck categories. The
main source of inspiration for our results was \cite{z}. However,
even for results that sound similarly, as for example our key
Theorem \ref{sumadirectaregulat}, stating that a finite direct sum
of $U$-regular objects in an abelian category is also $U$-regular,
and extending  \cite[Theorem 2.8]{z}, the proof is consistently
different. For some results we use Mitchell Theorem to reduce to
categories of modules, but most of the proofs are done inside the
general abelian category. We note that defining the concept of a
relative regular object in an arbitrary category presents the
advantage that we may pass to the dual category, and this is a
completely new method if we compare to the techniques used when
dealing with regular modules over rings. An example to illustrate
this statement is Corollary 2.8, which follows from Theorem 2.7 by
transferring the result to the dual category. As applications, we
give some results for the category of comodules over a coalgebra
and for the category of graded modules over a graded ring. We also
show that there is a link between the concept of a relative
regular object in a category and the theory of generalized
inverses of matrices.

We refer the reader to \cite{fr} and \cite{gab} for elements of
Category Theory, to \cite{dnr} for facts about coalgebras, and to
\cite{nv} for definitions and results about graded rings.

\section{Regular objects and basic properties}

Let $f:U\ra M$ and $h:M\ra U$ be morphisms in a category $\cal A$.
We say that $h$ is a generalized inverse of $f$ if $f=f\circ
h\circ f$ and $h=h\circ f\circ h$. It is easy to see that the
morphism $f:U\ra M$ has a generalized inverse if and only if there
exists a morphism $g:M\ra U$ such that $f=f\circ g\circ f$.
Indeed, an easy computation shows that $h=g\circ f\circ g$ is a
generalized inverse of $f$.

\begin{de}
Let $U$ and $M$ be objects of a category $\cal A$. Then $M$ is
called $U$-regular if and only if any morphism $f:U\ra M$ has a
generalized inverse.
\end{de}

\begin{re} \label{remarcaduala}
(1) If $U$ and $M$ are objects of a category $\cal A$, then $M$ is
$U$-regular in $\cal A$ if and only if $U$ is $M$-regular in the
dual category ${\cal A}^0$.\\
(2) If we associate to a ring $R$ a category with one object, such
that the elements of the ring are the morphisms and the
composition of morphisms is just the multiplication of $R$, then
$R$ is a von Neumann regular ring if and only if every morphism of
the associated category has a generalized inverse.
\end{re}

\begin{pr} \label{properties}
Assume that $M$ is $U$-regular in a category $\cal A$. The following assertions hold true.\\
(1) If $\pi :U\ra U''$ is an epimorphism in $\cal A$, then $M$ is also $U''$-regular.\\
(2) If $i:M'\ra M$ is a monomorphism in $\cal A$, then $M'$ is $U$-regular.\\
(3) If $\cal A$ is an additive category and $U'$ is a direct
summand of $U$, then $M$ is also $U'$-regular.
\end{pr}
{\bf Proof:} (1) Let $f:U''\ra M$ be a morphism. Since $M$ is
$U$-regular and $f\circ \pi:U\ra M$, there exists a morphism
$g:M\ra U$ such that $f\circ \pi =(f\circ \pi )\circ g\circ
(f\circ \pi )$.
Since $\pi$ is an epimorphism, we get that $f=f\circ (\pi \circ g)\circ f$.
Thus $M$ is $U''$-regular.\\
(2) Let $f:U\ra M'$ be a morphism. Since $M$ is $U$-regular and
$i\circ f:U\ra M$, there exists $g:M\ra U$ such that $i\circ
f=(i\circ f)\circ g\circ (i\circ f)$. Since $i$ is a monomorphism,
we see that $f=f\circ (g\circ i)\circ f$. Thus $M'$ is $U$-regular. \\
(3) Let $i:U'\ra U$ and $\pi :U\ra U'$ be morphisms such that $\pi
\circ i=1_{U'}$. If $f:U'\ra M$ is a morphism, then $f\circ \pi
:U\ra M$, so there exists $g:M\ra U$ such that $f\circ \pi=(f\circ
\pi )\circ g\circ (f\circ \pi )$. Then $f\circ \pi \circ i=(f\circ
\pi )\circ g\circ (f\circ \pi )\circ i$, showing that $f=f\circ
(\pi \circ g)\circ f$. Thus $M$ is $U'$-regular. \qed

\begin{pr} \label{transfer}
Let $F:{\cal A}\ra {\cal B}$ be a functor which is full and
faithful, and let $U,M$ be objects of $\cal A$. Then $M$ is
$U$-regular if and only if $F(M)$ is $F(U)$-regular.
\end{pr}
{\bf Proof:} Assume that $M$ is $U$-regular. Let
$\overline{f}:F(U)\ra F(M)$. Then there exists $f:U\ra M$ such
that $\overline{f}=F(f)$. Since $M$ is $U$-regular, there exists
$g:M\ra U$ such that $f=f\circ g\circ f$. Then
$\overline{f}=F(f)\circ F(g)\circ F(f)=\overline{f}\circ F(g)\circ \overline{f}$.
Thus $F(M)$ is $F(U)$-regular.\\

For the converse, assume that $F(M)$ is $F(U)$-regular and let
$f:U\ra M$. Then $F(f):F(U)\ra F(M)$, so there exists
$\overline{g}:F(M)\ra F(U)$ such that $F(f)=F(f)\circ
\overline{g}\circ F(f)$. Now there is $g:M\ra U$ such that
$\overline{g}=F(g)$, and then we see that $f=f\circ g\circ f$,
showing that $M$ is $U$-regular. \qed

\begin{ex}
Let $M$ and $U$ be sets. Then $M$ is $U$-regular in the category
$Set$ of sets. Indeed, let $f:U\ra M$ be a map. For any $m\in
Im(f)$ fix $a_m\in U$ such that $m=f(a_m)$. Also fix some element
$a\in U$, and define the map $g:M\ra U$ such that $g(m)=a_m$ for
any $m\in Im(f)$, and $g(m)=a$ for any $m\notin Im(f)$. Then
$f(g(f(u)))=f(a_{f(u)})=f(u)$ for any $u\in U$, so $f=f\circ
g\circ f$.\end{ex}

\begin{re} \label{remendreg}
(1) The definition of regular objects shows immediately that if
$U$ is an object of the additive category $\cal A$, then $U$ is
$U$-regular if and only if the
endomorphism ring $End(U)$ is a regular ring.\\
(2) If $M$ is $U$-regular, then any epimorphism $f:U\ra M$ splits,
i.e. there exists a morphism $g:M\ra U$ such that $f\circ g=Id_M$.
\end{re}

\section{Regular objects in abelian categories}

In this section we study regular objects in abelian categories. A
key characterization is the following.

\begin{pr}
Let $M$ and $U$ be objects of an abelian category $\cal A$. Then
$M$ is $U$-regular if and only if $Ker(f)$ is a direct summand of
$U$ and $Im(f)$ is a direct summand of $M$ for any morphism
$f:U\rightarrow M$.
\end{pr}
{\bf Proof:} Assume that $M$ is $U$-regular, and let $f:U\ra M$ be
a morphism. Then there exists a morphism $g:M\ra U$ such that
$f=f\circ g\circ f$. Let $f':U\ra Im(f)$ be the corestriction of
$f$ and $j:Im(f)\ra M$ be the inclusion morphism. Thus $f=j\circ
f'$. Hence $j\circ f'=j\circ f'\circ g\circ j\circ f'$. Since $j$
is a monomorphism we have that $f'=f'\circ g\circ j\circ f'$. But
$f'$ is an epimorphism, so $(f'\circ g)\circ j=Id_{Im(f)}$,
showing that $Im(f)$ is a direct summand in $M$. \\
On the other hand $f'\circ (g\circ j)=Id_{Im(f)}$ shows that
$Ker(f')=Ker(f)$ is a direct summand in $U$.\\

For the converse, to show that $M$ is $U$-regular, let $f:U\ra M$
be a morphism. As above let $f':U\ra Im(f)$ be the corestriction
of $f$ and $j:Im(f)\ra M$ be the inclusion morphism. Since
$Ker(f)$ is a direct summand of $U$, there exists $\beta :Im(f)\ra
U$ such that $f'\circ \beta =Id_{Im(f)}$. Also, since $Im(f)$ is a
direct summand in $M$, there exists $\al :M\ra Im(f)$ such that
$\al \circ j=Id_{Im(f)}$. Define $g:M\ra U$ by $g=\beta \al$. Then
\bea
f\circ g\circ f&=&f\circ \beta \circ \al \circ f\\
&=&j\circ f'\circ \beta \circ \al \circ j\circ f'\\
&=&j\circ f'\\
&=&f\eea We conclude that $M$ is $U$-regular. \qed

\begin{co}
Let $M$ and $U$ be objects of an abelian category $\cal A$. The
following assertions hold.\\
(1) If $U$ is injective, then $M$ is $U$-regular if and only if
$Im(f)$ is a direct summand of
$M$ for any morphism $f:U\rightarrow M$.\\
(2) If $M$ is projective, then $M$ is $U$-regular if and only if
$Ker(f)$ is a direct summand of $U$ for any morphism
$f:U\rightarrow M$.
\end{co}

\begin{co} \label{corsemisimple}
If $\cal A$ is a semisimple abelian category, then $M$ is
$U$-regular for any objects $M$ and $U$ of $\cal A$. Otherwise
stated, every morphism in a semisimple abelian category has a
generalized inverse.
\end{co}

\begin{re}
As the referee showed us, the concept of relative regular object
in a category can be related to the theory of generalized inverses
of matrices. Let $A$ be a $n\times m$-matrix and $B$ be a $m\times
n$-matrix. Then $A$ is called a generalized inverse of $B$ if
$ABA=A$ and $BAB=B$ (see for example \cite[Sec. 12.7]{lt}). \\
 If we apply Corollary
\ref{corsemisimple} to the category of finite dimensional vector
spaces over a field, we obtain the classical result that every
matrix has a generalized inverse (see \cite[p. 428, Prop. 1]{lt}).
\end{re}

\begin{co} \label{submodulsumand}
Let $N$ be a subobject of an object $M$ in an abelian category. If
$M$ is $N$-regular, then $N$ is a direct summand of $M$.
\end{co}

We are now interested to study finite direct sums of regular
objects in an abelian category. The following two results will be
of help for reducing the study to categories of modules.

\begin{pr} \label{result1}(\cite[Proposition 11.7]{fr})
Let ${\cal D}_0$ be a small subcategory of an abelian category
$\cal D$. Then there exists a small abelian full subcategory
${\cal D}'$ of $\cal D$ such that ${\cal D}_0$ is a subcategory of
${\cal D}'$.
\end{pr}

\begin{te} \label{result2}(Mitchell Theorem, see \cite[Theorem 11.6]{fr})
Let $\cal C$ be a small abelian category. Then there exist a ring
$A$ and a full and faithful exact functor $H:{\cal C}\ra
(Mod-A)^0$.
\end{te}

Now we can prove the main result of this section.

\begin{te} \label{sumadirectaregulat}
Let $M_1,M_2,\ldots, M_n$ be $U$-regular objects of an abelian
category $\cal A$, where $U$ is an object of $\cal A$. Then
$M_1\oplus M_2\oplus \ldots \oplus M_n$ is $U$-regular.
\end{te}
{\bf Proof:} Proposition \ref{result1}, Theorem \ref{result2} and
Proposition \ref{transfer} show that it is enough to prove the
result for the case where ${\cal A}=R-mod$, a category of modules
over a ring $R$. Using an inductive argument, it is enough to
prove for $n=2$.  Thus let $M_1$ and $M_2$ be $R$-modules which
are $U$-regular for a certain $R$-module $U$. Let $f:U\ra
M_1\oplus M_2$ be a morphism of $R$-modules. If $\pi _i:M_1\oplus
M_2\ra M_i$, $i=1,2$ are the natural projections, let $f_1=\pi
_1\circ f$ and $f_2=\pi _2\circ f$. Then $Ker(f)=Ker(f_1)\cap
Ker(f_2)$. Since $M_1$ is $U$-regular we have that $Ker(f_1)$ is a
direct summand in $U$. Since $M_2$ is $U$-regular, then by
Proposition \ref{properties} we see that $Ker(f_2\circ
i)=Ker(f_1)\cap Ker(f_2)$ is a direct summand of $Ker(f_1)$, where
$i:Ker(f_1)\ra U$ is the inclusion map. Therefore
$Ker(f)=Ker(f_1)\cap Ker(f_2)$ is a direct summand of $U$.\\

Now we show that $Im(f)$ is a direct summand in $M_1\oplus M_2$.
Clearly $Im(f)\subseteq Im(f_1)+Im(f_2)=Im(f_1)\oplus Im(f_2)$.
Since $Im(f_i)$ is a direct summand of $M_i$, $i=1,2$, we have
that $Im(f_1)\oplus Im(f_2)$ is a direct summand of $M_1\oplus
M_2$.
Thus it is enough to prove that $Im(f)$ is a direct summand of
$M_0=Im(f_1)\oplus Im(f_2)$. \\

The morphism $\pi _1$ induces the exact sequence

$$0\lra Ker(\pi_1')\stackrel{j}{\lra} Im(f)\stackrel{\pi_1'}{\lra} Im(f_1)\lra 0$$
where $\pi_1'$ is the restriction of $\pi_1$. We have the split
exact sequence

$$0\lra Ker(f_1)\lra U\stackrel{f_1}{\lra} Im(f_1)\lra 0$$
Then there exists $g_1:Im(f_1)\ra U$ such that $f_1\circ
g_1=Id_{Im(f_1)}$. Let $\al =f\circ g_1:Im(f_1)\ra Im(f)$. Then
$\pi_1'\circ \al =\pi_1'\circ f\circ g_1=f_1\circ
g_1=Id_{Im(f_1)}$,
so $Im(\al )$ is a direct summand of $Im(f)$. \\

If $\pi _1'':M_0\ra Im(f_1)$ is the restriction of $\pi_1$, then
clearly $Ker(\pi _1')\subseteq Ker(\pi_1'')=Im(f_2)$. Since $\pi
_1'\circ \al=Id_{Im(f_1)}$, there exists $j':Im(f)\ra Im(f_1)$
such that $j'\circ j=Id_{Ker(\pi _1')}$. Define $h=j_2\circ u\circ
j'\circ f:U\ra M_2$, where $j_2:Im(f_2)\ra M_2$ is the inclusion
map. Since $j'$ is an epimorphism we have that
$Im(h)=u(Ker(\pi_1'))=Ker(\pi_1')$. Now since $M_2$ is
$U$-regular, we have that $Im(h)$ is a direct summand of $M_2$. As
$Im(f_2)$ is a direct summand of $M_2$, we see that $Im(h)$
is a direct summand of $Im(f_2)$.\\

We have that $Im(f)=Ker(\pi_1')\oplus Im(\al )$. By the above
considerations we can define a morphism $\mu :Im(f_2)\ra Im(f_2)$
such that $\mu ^2=\mu$ and $Im(\mu )=Ker(\pi_1')$. Define $\theta
=i\circ \mu \circ \pi_2:M_0\ra M_0$, where $i:Im(f_2)\ra M_0$ is
the inclusion morphism. Then \bea
\theta ^2&=&i\circ \mu \circ \pi_2\circ i\circ \mu \circ \pi_2\\
&=&i\circ \mu^2\circ \pi_2\\
&=&i\circ \mu\circ \pi_2\\
&=&\theta \eea
Moreover we have that $Im(\theta)=Im(\mu )=Ker(\pi_1')$. \\

We also define $\gamma =i_0\circ \al \circ \pi_1:M_0\ra M_0$,
where $i_0:Im(f)\ra M_0$ is the inclusion map. Clearly $Im(\gamma
)=Im(\al )$. We have that \bea
\gamma ^2&=&i_0\circ \al \circ \pi_1\circ i_0\circ \al \circ \pi_1\\
&=&i_0\circ \al \circ \pi_1'\circ \al \circ \pi _1\\
&=&i_0\circ \al \circ \pi_1\\
&=&\gamma \eea We have that $Im(\gamma )\cap Im(\theta )=Im(\al
)\cap Ker(\pi_1')=0$. Moreover $\gamma \circ \theta=i_0\circ \al
\circ \pi_1\circ i\circ \mu\circ \pi_2=0$. By \cite[Lemma 1.1]{z}
we have that $Im(\gamma )+Im(\theta )=Im(\gamma )\oplus Im(\theta
)$ is a direct summand in $M_0$, and the proof is finished. \qed

By using Remark \ref{remarcaduala}(1) we get the following.

\begin{co} \label{sumadirectaduala}
Let $M,N_1,\ldots ,N_n$ be objects of an abelian category $\cal A$
such that $M$ is $N_i$-regular for any $1\leq i\leq n$. Then $M$
is $N_1\oplus \ldots \oplus N_n$-regular.
\end{co}

\begin{ex}
If we have an infinite family $(N_i)_{i\in I}$ of objects such
that $M$ is $N_i$-regular for any $i\in I$, then we do not
necessarily have that $M$ is $\oplus _{i\in I}N_i$-regular. To see
this, we consider \cite[Example 3.1]{z}, where it is given an
example of a von Neumann regular ring $R$, and a left ideal $J$ of
$R$ which is a regular $R$-module, while $End_R(J)$ is not a von
Neumann regular ring. Let $I$ be an infinite set such that there
is an epimorphism $R^{(I)}\ra J$. We claim that $J$ is not
$R^{(I)}$-regular. Indeed, if $J$ would be $R^{(I)}$-regular, then
by Proposition \ref{properties}(1) we have that $J$ is
$J$-regular. Hence by Remark \ref{remendreg}(1) we would get that
$End_R(J)$ is von Neumann regular, a contradiction.
\end{ex}

We remind that an abelian category has the property $(AB3)$ if it
has arbitrary direct sums.

\begin{co}
Let $U$ be a finitely generated object of an abelian category
$\cal A$ with $(AB3)$, and let $(M_i)_{i\in I}$ be a family of
objects of $\cal A$ such that $M_i$ is $U$-regular for any $i\in
I$. Then $\oplus_{i\in I}M_i$ is $U$-regular.
\end{co}
{\bf Proof:} Let $f:U\ra \oplus _{i\in I}M_i$ be a morphism. Then
there is a finite subset $J$ of $I$ such that $Im(f)$ is a
subobject of $\oplus _{i\in J}M_i$. Let $f_0:U\ra \oplus _{i\in
J}M_i$ be the corestriction of $f$. Since $\oplus _{i\in J}M_i$ is
$U$-regular, there exists $g_0:\oplus _{i\in J}M_i\ra U$ such that
$f_0=f_0\circ g_0\circ f_0$. Let $g:\oplus _{i\in I}M_i\ra U$
arising from the morphism $g_0$ and the zero morphism
$0:\oplus_{i\in I-J}M_i\ra U$. Then $f=f\circ g\circ f$, and this
ends the proof. \qed

We remind that in an  abelian category $\cal A$ with $(AB3)$, an
object $M$ is called $U$-generated if there exist a set $I$ and an
epimorphism $U^{(I)}\ra M$. If $I$ is finite (and $\cal A$ is just
an abelian category), then we say that $M$ is $U$-finitely
generated.

\begin{co}
Let $U$ be a projective object, and let $M$ be a $U$-regular
object of an abelian category.
The following assertions hold.\\
(i) If $M$ is $U$-finitely generated, then $M$ is projective.\\
(ii) If $\cal A$ has the property $(AB3)$ and $M$ is
$U$-generated, then $M$ is an inductive limit of projective
objects.
\end{co}
{\bf Proof:} (i) Let $f:U^n\ra M$ be an epimorphism. By Corollary
\ref{sumadirectaduala}, $M$ is $U^n$-regular. Hence by Remark
\ref{remendreg}(2), we see that $f$ splits, so then $M$ is a
direct summand of $U^n$. We conclude that $M$ is projective.\\
(ii) Since there exists an epimorphism $f:U^{(I)}\ra M$, we have
that $M$ is an inductive limit of subobjects that are epimorphic
images of objects of the form $U^n$ for some positive integers
$n$. Now the result follows from (i).\qed
\begin{lm}  \label{lemacentru}
Let $\cal A$ be an abelian category, $M$ an object of  ${\cal A}$,
and $\al$ an element of the center of $End_{\cal A}(M)$. Then
there exists $\beta$ in the center of $End_{\cal A}(M)$ with $\al
=\al\circ \beta\circ \al$ if and only if $M=Im(\al )\oplus Ker(\al
)$.
\end{lm}
{\bf Proof:} It follows from \cite[Lemma 3.3]{z} and Theorem
\ref{result2}. \qed

\section{Regular objects in locally finitely generated Grothendieck categories}

Throughout this section $\cal A$ is a Grothendieck category which
is locally finitely generated, i.e. it has a family $(U_i)_{i\in
I}$ of finitely generated generators.

\begin{de}
An object $M$ of $\cal A$ is called a regular object if $M$ is
$U_i$-regular for any $i\in I$.
\end{de}

To make the definition consistent, we need the following.

\begin{pr} \label{independent}
The concept of a regular object is independent on the choice of
the family of finitely generated generators of $\cal A$.
\end{pr}
{\bf Proof:} Assume that $M$ is $U_i$-regular for any $i\in I$.
Let $(V_j)_{j\in J}$ be another family of finitely generated
generators of $\cal A$. Let $j\in J$. Then there exist $i_1,\ldots
,i_n\in I$ and an exact sequence $U_{i_1}\oplus \ldots \oplus
U_{i_n}\longrightarrow V_j\longrightarrow 0$. By Corollary
\ref{sumadirectaduala} we have that $M$ is $U_{i_1}\oplus \ldots
\oplus U_{i_n}$-regular. Now by Proposition \ref{properties} (1)
we have that $M$ is $V_j$-regular. \qed

\begin{re} \label{remfgreg}
The proof of Proposition \ref{independent} shows that if $M$ is a
regular object of $\cal A$, then $M$ is $V$-regular for any
finitely generated object $V$ of $\cal A$.
\end{re}

\begin{ex}
If ${\cal A}=R-mod$, then a regular object in $\cal A$ is exactly
a regular $R$-module in the sense of \cite{z}.
\end{ex}

The following gives some properties of regular objects.

\begin{te}
Let $M$ be a regular object of a locally finitely generated
Grothendieck category $\cal A$. The following
assertions hold.\\
(1) The Jacobson radical $J(M)$ of $M$ is zero.\\
(2) The singular subobject $Z(M)$ of $M$ is zero.\\
(3) If $M$ is noetherian or artinian, then $M$ is semisimple.
\end{te}
{\bf Proof:} (1) Assume that $J(M)\neq 0$. Then there exists a
finitely generated non-zero subobject $N$ of $M$ such that
$N\subseteq J(M)$. By Corollary \ref{submodulsumand} we see that
$N$ is a direct summand of $M$, so there exists $P$ such that
$M=N\oplus P$. Since $N$ is finitely generated, there exists a
proper maximal subobject $N'$ of $N$. Then $N'\oplus P$ is maximal
in $M$, so
$J(M)\subseteq N'\oplus P$. Hence
$J(M)=J(M)\cap (N'\oplus P)=N'\oplus (J(M)\cap P)$. \\
On the other hand $J(M)=J(M)\cap (N\oplus P)=N\oplus (J(M)\cap
P)$, since $N\subseteq J(M)$ and the lattice of subobjects of $M$
is modular. Thus we have that $N'\oplus (J(M)\cap P)=N\oplus
(J(M)\cap P)$,
which shows that $N=N'$, a contradiction. Thus we must have $J(M)=0$.\\

(2) Assume that $Z(M)\neq 0$. Then there exists a non-zero
subobject $N$ of $M$ such that $N\simeq X/Y$ for some $X$ and $Y$
such that $Y$ is essential in $X$. Let $(U_i)_{i\in I}$ be a
family of finitely generated generators of $\cal A$. Then there is
$i\in I$ and a morphism $f:U_i\ra X$ such that $Im(f)$ is not a
subobject of $Y$. Let $\pi :X\ra X/Y$ be the natural projection,
and let $g=\pi \circ f:U_i\ra X/Y$. Clearly $Ker(g)=f^{-1}(Y)$.
Since $Y$ is essential in $X$, we have that $f^{-1}(Y)$ is
essential in $U_i$. Since $U_i$ is finitely generated, we have
that $Im(g)$ is also finitely generated. Now $Im(g)\subseteq
X/Y\simeq N$, so $Im(g)$ is a regular object. Hence $Ker(g)$ is a
direct summand of $U_i$. Since $Ker(g)$ is essential in $U_i$, it
follows that $Ker(g)=U_i$, so $g=0$, and then $Im(f)\subseteq Y$,
a contradiction. This shows that $Z(M)=0$.\\

(3) Let $M$ be noetherian. Then $M$ is finitely generated and
there exists a maximal subobject $N_1$ of $M$. Since $N_1$ is also
finitely generated, $M$ is $N_1$-regular and Corollary
\ref{submodulsumand} shows that $N_1$ is a direct summand of $M$.
Let $M=N_1\oplus S_1$, where $S_1$ is a simple object. In the same
way starting with $N_1$, we find that $N_1=N_2\oplus S_2$ for some
$N_2$ and a simple object $S_2$. We continue recurrently, and
since $M$ is noetherian, we end with $M=S_1\oplus S_2\oplus \ldots
\oplus S_n$
for some simple objects $S_1,S_2,\ldots ,S_n$. This shows that $M$ is semisimple.\\
If $M$ is artinian, let $s(M)$ be the socle of $M$. We have that
$s(M)$ is a finite direct sum of simple objects since $M$ is
artinian. Since $s(M)$ is also regular and essential in $M$, we
must have $s(M)=M$, showing that $M$ is semisimple. \qed

\begin{pr}
Let $\cal A$ be a locally finitely generated Grothendieck category.
The following assertions are equivalent.\\
(1) $\cal A$ is semisimple.\\
(2) Any object $M$ of $\cal A$ is regular.
\end{pr}
{\bf Proof:} Let $i\in I$ and $X$ be an essential subobject of
$U_i$. Then $U_i/X$ is $U_i$-regular, so $X=Ker(\pi )$ is a direct
summand of $U_i$, where $\pi :U_i\ra U_i/X$ is the natural
projection. Hence $X=U_i$, so the associated Goldie torsion theory
$\cal G$ is 0. But ${\cal A}/{\cal G}$ is a spectral category, so
${\cal A}$ is semisimple. \qed

\begin{te}
Let $\cal A$ be a locally finitely generated Grothendieck
category, and let $M$ be a regular object of $\cal A$. Then
$Z(End_{\cal A}(M))$ is a regular ring.
\end{te}
{\bf Proof:} Since $M$ is a regular object, $M$ is $N$-regular for
any finitely generated object $N$. Fix some $f$ in the center of
$End_{\cal A}(M)$. Let $M'$ be a finitely generated subobject of
$M$, and let $i:M'\ra M$ be the inclusion morphism. Since $M$ is
$M'$-regular and $f\circ i:M'\ra M$, there exists $g:M\ra M'$ such
that $f\circ i=(f\circ i)\circ g\circ (f\circ i)$. Since $f$ is in
the center of $End_{\cal A}(M)$, we have that $f\circ (i\circ
g)=(i\circ g)\circ f$, so then
\begin{equation} \label{form1}
f\circ i=f^2\circ (i\circ g\circ i)
\end{equation}
and also
\begin{equation} \label{form2}
f\circ i=(i\circ g)\circ f^2\circ i
\end{equation}
Using equation (\ref{form1}) we see that $f(M')=f(f(i\circ g\circ
i)(M'))$, so $M'\subseteq Im(f)+Ker(f)$.
To see this, note that in general $f(X)=f(Y)$ implies that
$X\subseteq Y+Ker(f)$ and $Y\subseteq X+Ker(f)$. \\

Since $M$ is the union of its finitely generated subobjects, we
get that $M\subseteq Im(f)+Ker(f)$, hence $M=Im(f)+Ker(f)$. We
show that $Im(f)\cap Ker(f)=0$. Indeed, if $K=Im(f)\cap Ker(f)$
were non-zero, then $L=f^{-1}(K)$ is non-zero. Since $\cal A$ is
locally finitely generated, there exists a non-zero subobject $N$
of $L$ such that $0\neq f(N)\subseteq K$. Since $f(K)=0$ we get
that $f^2(N)=0$, and then by applying equation (\ref{form2}) we
obtain that $f(N)=0$, a contradiction. Therefore  $Im(f)\cap
Ker(f)=0$ and $M=Im(f)\oplus Ker(f)$. The result follows now from
Lemma \ref{lemacentru}. \qed

More information about the endomorphism ring of a regular object
is given by the following result.

\begin{pr} \label{endsemiprime}
Let $\cal A$ be a locally finitely generated Grothendieck
category, and let $M$ be a regular object of $\cal A$.
The following assertions hold.\\
(1) $End(M)$ is a semiprime ring.\\
(2) If $M$ is finitely generated, then $End(M)$ is a regular ring.
\end{pr}
{\bf Proof:} (1) Let $f\in End(M)$, $f\neq 0$. Then there exists a
finitely generated subobject $M'$ of $M$ such that $f(M')\neq 0$.
By Remark \ref{remfgreg} we see that $M$ is $M'$-regular. Let $f'$
be the restriction of $f$ to $M'$. Then there exists $g':M\ra M'$
such that $f'=f'\circ g'\circ f'$. Let $g:M\ra M$, $g=i\circ g'$,
where $i:M'\ra M$ is the inclusion morphism. Then the restriction
of $f\circ g\circ f$ to $M'$ is $f'\circ g'\circ f'=f'\neq 0$.
Thus $fEnd(M)f\neq 0$, and we conclude
that $End(M)$ is a semiprime ring. \\
(2) By Remark \ref{remfgreg} we have that $M$ is $M$-regular, and
then $End(M)$ is a regular ring by Remark \ref{remendreg}. \qed

\section{Applications to coalgebras and graded rings}

If $C$ is a coalgebra over a field, we denote by ${\cal M}^C$ the
category of right $C$-comodules. This is a locally finite
Grothendieck category, see \cite{dnr}.
\begin{te}
Let $C$ be a coalgebra. Then $C$ is a regular object in the
category ${\cal M}^C$ if and only if $C$ is cosemisimple.
\end{te}
{\bf Proof:} If $C$ is cosemisimple then clearly $C$ is a regular
object in ${\cal M}^C$. Conversely, assume that $C$ is regular in
${\cal M}^C$. Let $M$ be a right comodule of finite dimension.
Then there exists a monomorphism $u:M\ra C^n$ for some positive
integer $n$. Since $C$ is regular, then $C$ is $M$-regular, hence
$C^n$ is $M$-regular. Hence there exists a morphism $v:C^n\ra M$
such that $u=u\circ v\circ u$. Since $u$ is a monomorphism we get
that $v\circ u=Id_M$, so $M$ is isomorphic to a direct summand of
$C^n$. This implies that $M$ is an injective object in ${\cal
M}^C$, and therefore the category ${\cal M}^C$ is semisimple, i.e.
$C$ is a cosemisimple coalgebra. \qed

Let $G$ be a group with identity element $e$, and let $R=\oplus
_{\sigma \in G}R_{\sigma}$ be a $G$-graded ring. Denote by $R-gr$
the category of graded left $R$-modules, which is a locally
finitely generated Grothendieck category, see \cite{nv}. $R$ is
called von Neumann gr-regular if and only if for any
$x_{\sigma}\in R_{\sigma}$ there exists $y\in R$ such that
$x_{\sigma}=x_{\sigma}yx_{\sigma}$ (clearly we can assume that
$y\in R_{{\sigma}^{-1}}$). It is clear that $R$ is von Neumann
gr-regular if and only if $Rx_{\sigma}$ is a direct summand in $R$
as an object of the category $R-gr$. If $M$ is an object of $R-gr$
and $\sigma \in G$, the $\sigma$-suspension $M(\sigma )$ of $M$ is
the graded $R$-module which is equal to $M$ as an $R$-module, and
whose grading is given by $M(\sigma )_{\tau}=M_{\tau \sigma}$ for
any $\tau \in G$.

\begin{te} \label{aplth2}
Let $R$ be a $G$-graded ring. Then the following assertions are equivalent.\\
(1) $R$ is von Neumann gr-regular.\\
(2) $R(\sigma )$ is $R$-regular in the category $R-gr$ for any
$\sigma \in G$.
\end{te}
{\bf Proof:} $(1)\Rightarrow (2)$ Let $f:R\ra R(\sigma )$ be a
morphism in $R-gr$, and let $x_{\sigma}=f(1)$. Then there exists
$y\in R_{{\sigma}^{-1}}$ such that
$x_{\sigma}=x_{\sigma}yx_{\sigma}$.
Then the map $g:R(\sigma )\ra R$ defined by $g(r)=ry$ is a morphism
in $R-gr$ and $f=f\circ g\circ f$.\\
$(2)\Rightarrow (1)$ Let $x_{\sigma}\in R_{\sigma}$. If $f:R\ra
R(\sigma )$ is defined by $f(r)=rx_{\sigma}$, then there exists a
morphism $g:R(\sigma )\ra R$ in $R-gr$ such that $f=f\circ g\circ
f$. Then if $y=g(1)$ we obtain that
$x_{\sigma}=x_{\sigma}yx_{\sigma}$. \qed

\begin{co}
Let $R_{\sigma}=\oplus _{\sigma \in G}R_{\sigma}$ be a graded ring.
Then the following assertions hold.\\
(1) If $R$ is von Neumann gr-regular then $R_{\sigma}$ is
$R_e$-regular for any $\sigma \in G$. In particular
$R$ is $R_e$-regular.\\
(2) If $R$ is strongly graded and $R_{\sigma}$ is $R_e$-regular
for any $\sigma \in G$, then $R$ is von Neumann gr-regular.
\end{co}
{\bf Proof:} (1) Let $f:R_e\ra R_{\sigma}$ be a morphism of
$R_e$-modules, and let $x_{\sigma}=f(1)$. Since $R$ is gr-regular,
there exists $y\in R_{{\sigma}^{-1}}$ such that
$x_{\sigma}=x_{\sigma}yx_{\sigma}$.
Define $g:R_{\sigma}\ra R_e$ by $g(r)=ry$ for any $r\in R_{\sigma}$.
Then clearly $f=f\circ g\circ f$.\\
(2) It follows by using Theorem \ref{aplth2}. \qed

\begin{co}\label{aplco1}
Let $G$ be a finite group and $R_{\sigma}=\oplus _{\sigma \in
G}R_{\sigma}$ be a graded ring. Then $R$ is von Neumann gr-regular
if and only if the smash product $R\#G$ is a von Neumann regular
ring.
\end{co}
{\bf Proof:} Assume that $R$ is von Neumann gr-regular. Let
$U=\oplus _{\sigma \in G}R(\sigma)$. Since $R$ is a finitely
generated object in $R-gr$, Theorem \ref{aplth2} shows that $U$ is
$R$-regular in $R-gr$. Since the $\sigma$-suspension functor
$T_{\sigma}:R-gr\ra R-gr$, $T_{\sigma}(M)=M(\sigma )$ for any
$M\in R-gr$, is an isomorphism of categories, we get that
$U=U(\sigma )$ is $R(\sigma)$-regular for any $\sigma \in G$.
Since $G$ is finite, we see that $U$ is $\oplus _{\sigma \in
G}R(\sigma )$-regular, i.e. $U$ is $U$-regular. Hence
$End_{R-gr}(U)$ is a von Neumann regular ring. But
$End_{R-gr}(U)\simeq R\#G$, so $R\#G$ is also von Neumann regular.
The converse is straightforward. \qed

\begin{co}\label{aplco2}
Let $R_{\sigma}=\oplus _{\sigma \in G}R_{\sigma}$ be a ring graded
by the arbitrary group $G$. If $R$ is von Neumann gr-regular, then
$End_{R_e}(R)$ and $End_{R-gr}(U)$ are semiprime rings.
\end{co}
{\bf Proof:} By Corollary \ref{aplco1}, $R$ is $R_e$-regular, and
by  Corollary \ref{aplco2} we have that $U$ is a regular object in
$R-gr$. Now the result follows from Proposition
\ref{endsemiprime}.\qed

{\bf Acknowledgement.} We would like to thank the referee for
indicating us the connection to the theory of generalized inverses
of matrices and for some other remarks that improved the
presentation of the paper.


\begin{thebibliography}{99}


\bibitem{dnr}
S. D\u{a}sc\u{a}lescu, C. N\u{a}st\u{a}sescu, \c{S}. Raianu, Hopf
algebras: An introduction, Pure and Applied Math. {\bf 235}
(2000), Marcel Dekker.



\bibitem{fr}
P. Freyd, Abelian categories. An introduction to the theory of
functors, Harper's Series in Modern Mathematics Harper \& Row,
Publishers, New York 1964.

\bibitem{gab}
P. Gabriel, Des cat\'{e}gories abeliennes, Bull. Soc. Math. France
90 (1962), 323-448.

\bibitem{go}
K. R. Goodearl, Von Neumann regular rings, 2nd edition, Fl:
Krieger Publishing Company. xvi, (1991).

\bibitem{lt}
P. Lancaster and M. Tismensky, The theory of matrices, second
edition, Academic Press, Orlando, 1985.

\bibitem{nv}
C.N\u ast\u asescu and F. Van Oystaeyen, Graded Ring Theory, North
Holland, 1982.

\bibitem{z}
J. Zelmanowitz, Regular modules, Trans. Amer. Math. Soc. {\bf 163}
(1972), 341-355.

\end{thebibliography}
\end{document}